\documentclass[12pt, reqno]{amsart}
\usepackage{amsmath, amstext, amsbsy, amssymb, amscd}

\newtheorem{theorem}{Theorem}[section]
\newtheorem{lemma}[theorem]{Lemma}

\newtheorem{corollary}[theorem]{Corollary}
\theoremstyle{definition}
\newtheorem{definition}[theorem]{Definition}

\theoremstyle{remark}
\newtheorem{remark}[theorem]{Remark}

\numberwithin{equation}{section}

\newcommand{\Supp}{{\rm Supp}}

{\vskip-\lastskip\medskip
  \noindent
  {\em #1.}\enspace
  }%
{\qed\par\medskip
  }

\begin{document}
\title[On the $D$-dimension of a certain type of threefolds]
              {On the $D$-dimension of a certain type of threefolds}

\author[Jing  Zhang]{Jing  Zhang}
\address{Department of Mathematics, University of 
Missouri, Columbia, MO
65211, USA}
\email{zhangj@math.missouri.edu}

\subjclass{Primary  14J30; Secondary  32Q28}

\begin{abstract}
Let $Y$ be an algebraic manifold of dimension 3 with 
$H^i(Y, \Omega^j_Y)=0$
for all $j\geq 0$, $i>0$ and 
$h^0(Y, {\mathcal{O}}_Y) > 1$. 
Let $X$ be a smooth completion of  $Y$  such that
the boundary $X-Y$  is the support  of  an  effective divisor 
$D$   on $X$  with    simple  normal  crossings.  
We prove that the $D$-dimension of $X$ cannot be 2, i.e., 
either  any two nonconstant regular   functions   
are algebraically  dependent 
or there are three  algebraically  independent nonconstant
 regular functions on 
$Y$. 
Secondly, if the  $D$-dimension of $X$
is greater than 1, then the associated scheme of $Y$ is isomorphic to 
Spec$\Gamma(Y, {\mathcal{O}}_Y)$. 
Furthermore, we prove that  an   algebraic
manifold  $Y$   of  any dimension  $d\geq 1$ is affine if and only if
$H^i(Y, \Omega^j_Y)=0$
for all $j\geq 0$, $i>0$ and 
 it is regularly 
separable, i.e., for any two distinct points 
$y_1$,  $y_2$  on  $Y$, there is a regular function $f$ on $Y$ such that
$f(y_1)\neq f(y_2)$. 
\end{abstract}

\maketitle
\date{}
\section{Introduction}

We work over  complex  number  field  $\Bbb{C}$.

 We continue our investigation of  
 three  dimensional  algebraic 
 manifolds $Y$ with 
$H^i(Y, \Omega^j_Y)=0$   
     for all $j\geq 0$ and $i>0$, where $\Omega^j_Y$
is the sheaf of regular $j$-forms 
on $Y$. 
Originally   this  is a question   raised by J.-P. Serre for 
complex manifold  \cite{Se}, there $\Omega^j_Y$
is the sheaf of   holomorphic   $j$-forms 
on $Y$.
We are interested in the classification of  threefolds 
with this vanishing property.
  Let $X$ 
 be a smooth completion of  $Y$  such that
the boundary $X-Y$  is the support  of  an  effective divisor 
$D$   on $X$  with    simple  normal  crossings. 
Our approach   depends on the  the $D$-dimension 
of $X$,   
 a notion due to Iitaka  \cite{I1}. 
If for all $m>0$, $H^0(X, {\mathcal{O}}_X(mD))=0$, then 
the $D$-dimension
$\kappa(D, X)=-\infty$. Otherwise, 
$$\kappa(D, X)=tr.deg_{\Bbb{C}}\oplus_{m\geq 0}
H^0(X, {\mathcal{O}}_X(mD))-1.$$
The Kodaira dimension $\kappa(X)$ of $X$  is defined to be 
$\kappa(K_X, X)$, where $K_X$ is the conanical divisor 
of $X$.  
An equivalent definition of $D$-dimension 
and some properties are reviewed in Section 2.

In our previous papers \cite{Zh1, Zh2}, 
we understand  the structure of $Y$ very well when
$\kappa(D, X)=1$. We know  
if  $H^i(Y, \Omega^j_Y)=0$ for all $j\geq 0$ and $i>0$, then 
 $X-Y$
is connected. If we  also assume 
$\kappa(D, X)\geq  1$, then  
$Y$ contains no complete curves  and we have the following results.

    (1) There is a smooth projective 
curve $\bar{C}$,  and a smooth  affine curve $C$ such that 
the following diagram commutes

\[
  \begin{array}{ccc}
    Y                           &
     {\hookrightarrow} &
    X                                 \\
    \Big\downarrow\vcenter{%
        \rlap{$\scriptstyle{f|_Y}$}}              &  &
    \Big\downarrow\vcenter{%
       \rlap{$\scriptstyle{f}$}}      \\
C        & \hookrightarrow &
\bar{C},
\end{array}
\]
where $f$ is proper and surjective, every fibre of $f$ over $\bar{C}$
 is connected, and a 
 general fibre  is  smooth. Also a general fibre 
 of $f|_Y$ is connected  and smooth. 
In particular, 
every fibre $S$ of $f|_Y$ over $C$  satisfies the same vanishing 
       condition, i.e., 
$H^i(S, \Omega^j_S)=0$. In fact, all smooth fibres are of the same type.

(2) The Kodaira dimension $\kappa (X)=-\infty$ and the $D$-dimension
$\kappa(D, X)=1$  if a general fibre in the  above  fibre space 
$f|_Y: Y\rightarrow  C$  is not 
affine.\\

 In \cite{Zh2}, we prove that
     there exist nontrivial (i.e., nonaffine and nonproduct)
 threefolds $Y$ with 
$H^i(Y, \Omega^j_Y)=0$ for all $j\geq 0$ and $i>0$
such that    $\kappa (X)=-\infty$ and 
$\kappa (D, X)=1$.

 In this paper,
 we   will   show that  the $D$-dimension of $X$ cannot be 2. 
This means that 
if $h^0(Y, {\mathcal{O}}_Y) > 1$ (which is equivalent to 
$\kappa (D, X)\geq 1$), then
either any  two 
nonconstant regular   functions   on $Y$
are algebraically  dependent 
or there are three  algebraically  independent nonconstant
 regular functions on 
$Y$. 
Secondly, if the  $D$-dimension of $X$
is greater than 1, then $Y$ 
(we identify $Y$ with its associated scheme, \cite{H1},
Chapter II, Proposition 2.6) 
is   isomorphic   to Spec$\Gamma(Y, {\mathcal{O}}_Y)$. 
Furthermore,   
$Y$ is affine if and only if
 $Y$ is regularly separable. 
 
\definition{An algebraic variety $Y$ is regularly separable if 
for any two distinct points $y_1$ and $y_2$ on 
$Y$, there is a regular
function $f$ on $Y$ such that $f(y_1)\neq f(y_2)$.\\

 \noindent {\bf{Main Theorem }}
{\it  Let $Y$ be a smooth threefold with  
$H^i(Y, \Omega^j_Y)=0$ for all $j\geq 0$ and $i>0$.
Let $X$ be a smooth completion 
of $Y$ and $D$  be the effective boundary divisor
with simple normal crossings such that the support of
$D$ is $X-Y$. Then we have

(1)  $\kappa (D, X)\neq 2$.

 (2) If 
$\kappa (D, X)>1$, then the associated scheme of  $Y$ is isomorphic to 
\mbox{Spec}$\Gamma(Y, {\mathcal{O}}_Y)$. 

(3) $Y$  is  affine if and  only  if   $Y$  is  regularly
separable.
  }\\

With the same notion as in the Main Theorem, we have 
\begin{corollary} If $Y$ is a 
 smooth threefold with  
$H^i(Y, \Omega^j_Y)=0$ for all $j\geq 0$ and $i>0$, then the following
four conditions are equivalent

(1) $\kappa (D, X)>1$;

(2)  $\kappa (D, X)=3$;

(3)   $Y$  is  regularly
separable;

(4) $Y$  is  affine.   
\end{corollary}

The statement (3) in the Main Theorem 
holds for any dimensional algebraic  manifolds.

\begin{theorem} An algebraic manifold  $Y$   is    affine   if and only if 
 $H^i(Y, \Omega^j_Y)=0$
for all $j\geq 0$, $i>0$ and  $Y$ is 
regularly separable. 
 \end{theorem}

In order to prove the  Main Theorem, 
We need the following  result for surfaces,   which 
is interesting on its own.

\begin{theorem}
Let  $Y$  be an irreducible  quasi-projective   surface.
 Let $X$ be a projective  surface containing $Y$.
 Then $Y$ is affine if and only if
the following three conditions hold

(1) Y contains no complete curves;

(2)  The boundary $X-Y$ is connected;

(3) $\kappa (D, X)=2$, where $D$ is an effective divisor 
with support $X-Y$.   
\end{theorem}

In Section 2, we will    prove   some results   including Theorem 1.4
 for  surfaces
which will be used in Section 3. We will prove
the Main  Theorem and Theorem 1.3 in Section
3. 
The idea to   show  $\kappa (D, X)\neq 2$ is to
compute the $D$-dimension by the fibre space we constructed 
in   \cite{Zh1}  and a result of Fujita \cite{Fuj2}. 
To prove that $Y$ is 
 isomorphic to 
Spec$\Gamma(Y, {\mathcal{O}}_Y)$,
we    show   that there is an 
injective   birational  morphism 
 from
 $Y$ to 
Spec$\Gamma(Y, {\mathcal{O}}_Y)$. 
Then by  Zariski's Main Theorem (\cite{Mu1}, Chapter
3, Section 9)  and  a theorem of Neeman [N], $Y$ is affine. 
\\

\noindent 
{\bf{Acknowledgments}}  \quad  I would like to    thank  
   the following 
      professors for  helpful discussions:  Steven Dale Cutkosky,
     Dan Edidin,   N.Mohan  Kumar, 
      Zhenbo Qin,  and Qi Zhang.

\section{Surfaces}

We start with the definition of 
$D$-dimension and its properties we will use later. 
For general references  see (\cite{I4};  \cite{I1};   \cite{I2},  Leture 3;
 \cite{Uen}, Chapter 2).
Let $X$ be a normal  projective  variety and 
$D$ be  a Cartier
      divisor on $X$.  Then   associated to $D$ 
     we have  a line bundle
     ${\mathcal{O}}_X(D)$. If for all integers  $m> 0$ we have
     $H^0(X, {\mathcal{O}}_X(mD))=0$, then we define 
     the  $D$-dimension of $X$, denoted by $\kappa (D, X)$, to be $-\infty$.
     If  $h^0(X, {\mathcal{O}}_X(mD))\geq 1$ for some $m$, 
     choose a basis $\{f_0, f_1, \cdot \cdot\cdot, f_n\}$
     of the linear space 
     $H^0(X, {\mathcal{O}}_X(mD))$, it defines a rational 
     map 
     $\Phi _{|mD|}$
     from $X$ to the projective space 
     ${\Bbb{P}}^n$ by sending a point $x$ on $X$ to
     $(f_0(x), f_1(x), \cdot \cdot\cdot, f_n(x))$ in ${\Bbb{P}}^n$.  
     Then we define
     $\kappa (D, X)$ to be the maximal dimension of the images
      of the rational map  $\Phi _{|mD|}$, i.e., 
      \begin{equation}  \kappa (D, X)= \max_m\{\mbox {dim} (\Phi _{|mD|}(X))\}. 
      \end{equation}
If $X$ is not normal, let $\pi: X^*\rightarrow X$ be the normalization 
of $X$, then we define
\begin{equation}\kappa(D, X)=\kappa(\pi^*D, X^*),
 \end{equation}
where  $\pi^*D$
denotes the pull back of the  Cartier divisor $D$  
from $X$ to  $X^*$. 
The 
       $D$-dimension of a variety is  a birational invariant. 
        We   do not change the $D$-dimension
by blowing up or blowing  down (\cite{I4}, Section 5; 
\cite{Uen}, Chapter 2,
Theorem 5.13). 
        More precisely,  let  $f$: $X'\rightarrow X$
        be a surjective morphism between  two  complete
        varieties, let $D$ be a divisor on 
        $X$ and $E$ an effective divisor on $X'$ such that
        codim$f(E)\geq 2$, then  
        \begin{equation}\kappa (f^{-1}D+E, X')=\kappa (D, X), 
        \end{equation}
where  $f^{-1}D$   is the   reduced transform of $D$:
$f^{-1}D=\sum \Gamma_i$, $\Gamma_i$'s are the irreducible 
components  of $D$.  
Another property of $D$-dimension is that 
       it does not depend on the coefficients of $D$ 
       under  a mild condition which is certainly true in our case
       (since we always choose effective boundary divisor $D$
       with simple  normal crossings).
         Let $D_1$, $D_2$, $\cdot$$\cdot$$\cdot$, $D_n$
       be  any divisor on $X$ such that for every $i$, $0\leq i \leq n$, 
       $\kappa (D_i, X)\geq 0$, then for   integers 
       $p_1> 0,\cdot\cdot\cdot$, $p_n>0$,    we have
(\cite{I4}, Section 5)
       \begin{equation}\kappa(D_1+\cdot\cdot\cdot+D_n,X)=
       \kappa(p_1D_1+\cdot\cdot\cdot+p_nD_n,X).
\end{equation}       
With these two properties, 
we  may change the coefficients of $D$ to different positive integers 
or blow  up 
a curve or point on the boundary $X-Y$ freely and still call it $D$.

  A fibre space is a morphism $f$: $X\rightarrow  Z$ 
      which  is proper and surjective  with general fibre  connected.
       Suppose both $X$ and $Z$ are  nonsingular, then by 
       Theorem 5.11,  \cite{Uen},
       for any Cartier divisor $D$ on $X$,  there  exists an open 
       dense subset $U$ of $Z$  in complex topology such that for any fibre 
       $X_z=f^{-1}(z)$, $z\in U$, the inequality
       \begin{equation} \kappa (D, X)\leq \kappa (D_z, X_z)+\mbox {dim}(Z)
       \end{equation} 
       holds, where $D_z=D|_{X_z}$, the restriction divisor on 
the fibre $X_z$.  
\begin{lemma} Let  $Y$  be an irreducible smooth surface
without complete curves. Let $X$ be its smooth completion.
Suppose that the boundary $X-Y$ is connected. Let $D$ be an
effective divisor with support $X-Y$. If 
$\kappa (D, X)=2$, then  $X-Y$ is  the  support of an ample divisor $P$.
 \end{lemma}
{\it Proof}. Since $Y$ contains no complete curves,
the  boundary  $X-Y$  cannot  be blown down to a point. 
And $X-Y$  is  of pure codimension 1 since it is connected. 

Write the Zariski decomposition  $D=P+N$, where 
$N$ is negative  definite, $P$  is effective and nef and 
any prime component of  
$N$  does not intersect $P$  \cite{Za}. We may assume that 
both $P$ and $N$ are integral by multiplying a positive 
integer  to the equation since both $P$ and  $N$ are $\Bbb{Q}$ divisors
($D$ is a Weil divisor but $P$ and $N$  have  
rational coefficients). 
Let  $\Supp{D}=\{D_1, D_2, \cdot\cdot\cdot, D_n \}=X-Y$.
Since   $\kappa (D, X)=2$,  $P^2>0$ 
(\cite{Sa1}; \cite {Ba}, Corollary 14.18). 
First we claim that $\Supp{P}=\Supp{D}=X-Y$. If 
$\Supp{P}\neq X-Y$, then there is a prime component, say $D_1$, 
in $X-Y$ such that  $P\cdot D_1>0$ and  $D_1$ is not a component of 
$P$  since $X-Y$ is connected.   Let 
$$Q=mP+D_1,$$
 where $m$ is  a  big  positive integer.
Then  $Q$ is an effective divisor and $\Supp{Q}=\Supp{P}\cup D_1$. 
Since  $P^2>0$,  we may choose $m$ such  that  
$$ Q^2=m^2P^2+2mP\cdot  D_1+ D_1^2>0.  
$$
For every prime  component  $E$ in $P$, since $P$ is nef and 
$D_1$ is not contained in $\Supp {P}$, for  sufficiently large $m$,
 we have 
$$ Q\cdot E=mP\cdot  E+D_1\cdot E\geq 0, \quad \quad 
D_1\cdot  Q=mD_1\cdot P +D_1^2 >0.
$$
Since $Y$  contains not complete curves, 
any irreducible complete curve outside $X-Y$ intersects 
$X-Y$.  
Thus we get a new effective divisor $Q$ 
such  that  $Q$ is nef and $Q^2>0$. We may replace $P$ by $Q$
and still call it $P$. 
By finitely many such   replacements, we 
can   find an effective nef  divisor $P$     such that 
$P^2>0$ and $\Supp{P}=\Supp{D}=X-Y$.

We  claim  that the boundary $X-Y$ is the support of an ample divisor. 
In fact, the following three conditions imply  the  ampleness: 

(1) $X-Y$ is connected;

(2) $Y$ contains no complete curves;

(3) There is an effective nef divisor $P$ with supp$P=X-Y$ and $P^2>0$.

If $P$ is not ample, then there is an irreducible complete 
curve $C$ in $X$ such that
$P\cdot C=0$ by Nakai-Moizshon's ampleness criterion
(\cite{H1}, Chapter V). Since $Y$ has no complete curves, 
$C$ must be one of the $D_i's$. Rearrange the order, we may assume 
$D_i\cdot P=0$ for all $i$,  $i=1,2,..., r$ and $D_j\cdot P>0$
for  all  $j$,  $j= r+1,..., n$.  
Write 
$$P=\sum_{i=1}^ra_iD_i+\sum_{j=r+1}^nb_jD_j=A+B,$$
where
$A=\sum_{i=1}^ra_iD_i$, $B=\sum_{j=r+1}^nb_jD_j$. 
Then for all $i$,  $i=1,...r$,
$$0=P\cdot D_i = A\cdot  D_i+B\cdot D_i.$$
Since $D_i$ is not a component of $B$,  for all $i$,  $i=1,..., r$,
$B\cdot D_i\geq 0$. So $A\cdot  D_i\leq  0$  for every 
$i=1,..., r$. Thus the intersection matrix 
$[D_s\cdot D_t]_{1\leq s, t\leq r}$
is  negative  semi-definite \cite{Art}.
Since $A\cup B=X-Y$ is connected, there is at least one 
component $D_{i_o}$  of $A$,  such that 
$D_{i_o}\cdot B>0$. Hence 
 $D_{i_o}\cdot A<0$. This implies that 
the intersection matrix 
$[D_s\cdot D_t]_{1\leq s, t\leq r}$
is  negative  definite \cite{Art}.
Therefore there is an effective divisor 
$E=\sum_{i=1}^r\alpha _iD_i$  such that 
$E\cdot D_i<0 $ for all $i$, 
$i=1,..., r$ \cite{Art}.

So there are positive numbers 
$\alpha_i$, $i=1,..., r$ such that for every $i$, 
$1\leq i\leq r$, $E\cdot D_i< 0$,
where $E= \sum_{i=1}^r\alpha_iD_i$.  Let $P_1=mP-E$, $m\gg 0$, then
$P_1^2>0$, $P_1$ is nef  and if $1\leq i\leq r$, 
$$P_1\cdot  D_i=-E\cdot D_i>0.$$
If $r+1 \leq j\leq n$, then  for sufficiently large $m$, 
$$  P_1\cdot  D_j=mP\cdot D_j-E\cdot D_j>0.
$$
Thus $P_1$ is an effective ample divisor with support $X-Y$. Replace 
$P$ by $P_1$, we  have   shown    that $X-Y$ is 
the support of an ample divisor 
$P$.  
\begin{flushright}
 Q.E.D. 
\end{flushright}

\begin{theorem}
Let  $Y$  be an irreducible open smooth surface.
 Let $X$ be its smooth completion. Then $Y$ is affine if and only if
the following three conditions hold

(1) Y contains no complete curves;

(2)  The boundary $X-Y$ is connected;

(3) $\kappa (D, X)=2$, where $D$ is an effective divisor 
with support $X-Y$. 
\end{theorem} 
{\it Proof}.  If $Y$ is affine, then the above three conditions hold 
 (\cite{Ba}, Corollary 14.18;   \cite{H2}, Chapter 2, Section 3 and 4).
 Conversely, if these three conditions hold, then by 
Lemma 2.1,  the boundary $X-Y$ is the support of an ample divisor $P$.  
Thus $Y$ is affine by a theorem of Goodman (\cite{H2}, Chapter 2, 
Theorem 4.2).
\begin{flushright}
 Q.E.D. 
\end{flushright}

Notice that the above theorem holds for complete normal surfaces. 
For a complete  normal surface $X$, the intersection theory  is due to
Mumford \cite{Mu2}.
Let Div$X$  be the group of  Weil divisors
of $X$.  Let  Div$(X, {\Bbb{Q}})=$Div$(X)\otimes  {\Bbb{Q}}$
be the  group of  $\Bbb{Q}$-divisors. 
The intersection pairing     
$${\mbox{Div}}(X, {\Bbb{Q}})\times {\mbox{Div}}(X, {\Bbb{Q}})
\rightarrow  \Bbb{Q}
$$
is defined  in the following way.  Let $\pi: X'\rightarrow X$
be a resolution and let $A=\cup E_i$
denote the exceptional  set of $\pi$.
For a $\Bbb{Q}$-divisor $D$ on $X$
we define the inverse image 
$\pi^*D$ as
$$ \pi^*D=\bar{D}+\sum a_iE_i 
$$
where $\bar{D}$  is the strict transform  of $D$ by $\pi$
and the rational numbers $a_i$ are uniquely determined by the equations 
$\bar{D}E_j+\sum a_iE_iE_j=0$  for all $j$. For two divisors 
$D$ and $D'$ on $X$,  define their intersection number  
$$   D\cdot D'=\pi^*D\cdot \pi^*D'.
$$

\begin{lemma}{\bf[Fujita]} Let  $D$ be an effective  $\Bbb{Q}$-divisor
on a normal projective  surface $X$. Then  there exists a unique 
decomposition 
$$D=P+N$$
satisfying the following conditions:

(1) N is an effective $\Bbb{Q}$-divisor and either N=0 or the intersection 
matrix of the irreducible components of N is negative definite;

(2) P is a nef $\Bbb{Q}$-divisor and the intersection of P 
with each irreducible 
component of N is zero. 
\end{lemma}

\begin{theorem}
Let  $Y$  be an irreducible  quasi-projective   surface.
 Let $X$ be a projective  surface containing $Y$.
 Then $Y$ is affine if and only if
the following three conditions hold

(1) Y contains no complete curves;

(2)  The boundary $X-Y$ is connected;

(3) $\kappa (D, X)=2$, where $D$ is an effective divisor 
with support $X-Y$.   
\end{theorem} 
{\it Proof}. If we have a surjective and finite morphism from 
a variety  $Y'$ to $Y$, then $Y$ is affine if and only if $Y'$ is affine 
by Chevalley's theorem (\cite{H2}, Chapter 2, Corollary 1.5). 
Thus $Y$ is affine if and only if its normalization is affine.
 So we may assume that both 
$Y$ and $X$ are normal by taking their normalization. On a normal 
 projective surface, the intersection theory and Zariski decomposition
remain true by  Lemma 2.3 \cite{Mu2,  Sa2}. Therefore Lemma 2.1 
holds for normal   projective 
surfaces.
In fact, write the  Zariski decomposition
$D=P+N$ as in the above Lemma 2.3, then 
$P^2>0$ (\cite{Ba}, Corollary  14.18, Page 222).
By the same argument 
as in the proof of  Lemma 2.1, we can find a 
new effective 
nef divisor, still denoted by $P$, such that  
supp$P=X-Y$.   By changing the coefficients of $P$,
we can  find     an ample divisor supported 
in  $X-Y$.   
\begin{flushright}
 Q.E.D. 
\end{flushright}
 
If $Y$ is an irreducible  smooth surface with 
$H^i(Y, \Omega^j_Y)=0$
for all $j\geq 0$ and $i>0$, then the first two conditions are 
satisfied \cite{Ku}. And by  Theorem 3.1 (next section), 
$\kappa(D, X)=0$ or $2$. 
 If $\kappa(D, X)=2$, then
 $Y$ is affine.

\begin{remark} The above Theorem 2.4 does not hold  for threefolds. 
If $Y$ is a smooth quasi-projective threefold without complete curves
and the boundary $X-Y$ is connected for a smooth completion $X$ of $Y$,
then the 
boundary   being  of pure codimension 1 and 
$\kappa(D, X)=3$ cannot guarantee that $Y$ is affine. The reason is that 
in surface case,
the two conditions, i.e., $Y$ contains no complete curves and
$X-Y$ is connected, imply that for any smooth completion $Z$ 
(may be  different from $X$)
of $Y$,
 the boundary $Z-Y$  is 
of pure codimension 1. This is of course  not true in higher dimension.
For instance, remove a hyperplane section $H$ and a line
$L$ from ${\Bbb{P}}^3$, where $L$ is not contained in $H$. 
Let $Y={\Bbb{P}}^3-H-L$.  Then $Y$ contains no complete curves.
 Let 
$f:X\rightarrow {\Bbb{P}}^3$
be the blowing up of ${\Bbb{P}}^3$ along $L$. Then $X$ is a  smooth 
projective threefold and  $Y$ is an open subset of $X$. 
Let $D=f^{-1}(H)+E$, where $E$ is the exceptional divisor.
Then by equation (2.3),
$\kappa (D, X)=\kappa (H, {\Bbb{P}}^3)=3$. But $Y$ is not affine 
since  the boundary ${\Bbb{P}}^3-Y$ is not of pure codimension
1  (\cite{H2}, Chapter 2, Proposition 3.1).  
 \end{remark}

For a projective manifold $M$, let $L$
be a line bundle on $M$, then  there is a  Cartier divisor 
$D$ determined by $L$. We define $\kappa(L,M)=\kappa(D,M)$.

\begin{lemma}[{\bf Fujita}] Let 
$M$ and  $S$  be two projective manifolds.  
 Let $\pi: M\rightarrow S$
  be a  fibre space and let $L$ and $H$ be  line bundles  on $M$
  and $S$ respectively. Suppose that  
  $\kappa (H, S)=\dim S$ and that
  $\kappa (aL-b\pi^*(H))\geq 0$
for certain  positive integers $a$, $b$.  Then
$\kappa(L,M)=\kappa (L|_F, F)+ \kappa (H, S)$
for a general  fibre $F$  of  $\pi$.   
\end{lemma}

\begin{corollary} Suppose that we have a surjective morphism from an 
irreducible  smooth quasi-projective surface $Y$ to a smooth affine curve $C$.
Let $X$ be a  smooth projective surface containing $Y$.
If      $Y$ contains no complete curves and  
the boundary $X-Y$ is connected, then 
$Y$ is affine.
 \end{corollary}
{\it Proof}. Let $f:Y\rightarrow C$ be the given morphism.
 Then  $f$ gives a rational map from $X$ to $\bar{C}$, where 
$\bar{C}$  is the smooth completion of $C$.  Resolve the
indeterminacy of $f$ on the boundary $X-Y$.
We may replace $X$ by its suitable blowing up and assume that
$f:X\rightarrow  \bar{C}$  is surjective and proper morphism. 
Notice that this  procedure  does not  change $Y$.
$Y$ is still an open subset of $X$. 
By Stein factorization, we may assume that every fibre is 
connected and general fibre is smooth.
Pick a point $t_1\in \bar{C}-C$, then 
$$h^1(\bar{C}, {\mathcal{O}}_{\bar{C}}(nt_1))
=0$$
since $nt_1$ is ample for large $n$
(\cite{H1}, Chapter IV, Corollary  3.3). 
By the  Riemann-Roch formula,
$$h^0(\bar{C}, {\mathcal{O}}_{\bar{C}}(nt_1))=
1+n-g(\bar{C}).$$
So 
$\kappa(t_1, \bar{C})=1$.  For a general point 
$t\in C$, by  Riemann-Roch, there is  a positive integer $m$,
 such that 
$h^0(\bar{C}, {\mathcal{O}}_{\bar{C}}(mt_1-t))>1$. 
 Let $s$ be a nonconstant section of 
$H^0(\bar{C}, {\mathcal{O}}_{\bar{C}}(mt_1-t))$, then
$$ {\mbox{div}} s +mt_1-t\geq 0.
$$
Pull it back to $X$, we have 
$$  f^*({\mbox{div}} s +mt_1-t)={\mbox{div}} f^*(s)+mf^*(t_1)-f^*(t)\geq 0.
$$
Let $D_1=f^*(t_1)$ and $F=f^*(t),$ then  
$h^0(X, {\mathcal{O}}(mD_1-F))>0$.  Choose an effective divisor 
$D$ with support $X-Y$ such that  $D_1\leq D$, then we  have
$$ h^0(X, {\mathcal{O}}(mD-F))\geq  h^0(X, {\mathcal{O}}(mD_1-F))>0. 
$$
Since  $F|_Y$ is a smooth affine curve (\cite{H2}, Chapter 2, Proposition 4.1),
$$ h^0(F, {\mathcal{O}}_F(mD|_F))\geq  n+1-g(F).
$$
Therefore
 $\kappa(D|_F, F)=1$.
By Lemma 2.6 and Equation  (2.4), 
$$\kappa(D, X)=\kappa(mD, X)=\kappa(mD|_F, F)+\kappa(t_1, \bar{C})=2.$$
By Theorem 2.2, $Y$ is affine.
\begin{flushright}
 Q.E.D. 
\end{flushright}
\begin{definition}
A   complex   space  $Y$
 is Stein if and only if 
 $H^i(Y, G)=0$  for every analytic coherent sheaf $G$ on $X$
 and all positive integers $i$.
\end{definition} 
 If $Y$ is a  holomorphic variety, then 
 $Y$ is Stein if and only if  it is both holomorphically 
 convex and holomorphically separable (\cite{Gu}, Page 143). 
We say that  $Y$ is  
holomorphically 
 convex if for any discrete sequence $\{y_n\}\subset Y$,
 there is a holomorphic function $f$ on $Y$ such that
 the supremum of the set $\{|f(y_n)|\}$ is 
 $\infty$. $Y$ is  
holomorphically  separable if for every pair $x,y\in Y$,
$x\neq y$,  there is a holomorphic function $f$ on $Y$ such that
$f(x)\neq f(y)$.

\begin{remark} Since Theorem 2.2 is not true for threefolds,
Corollary 2.7 does not hold for threefolds. We have the following 
counter-example.

 Let $C$ be an elliptic curve (smooth and projective) 
and $E$ the unique nonsplit 
      extension of $\mathcal{O}$$_C$ by itself.  
      Let ${Z=\Bbb{P}}_C(E)$ and  $D$ be the canonical section,
then  $H^i(S, \Omega^j_S)=0$ for all $i>0$ and $j\geq 0$, where  $S=Z-D$
\cite{Ku}.
Let  $F$ be a smooth affine curve and $Y=S\times F$, then  
$H^i(Y, \Omega^j_Y)=0$ by   K$\mbox{\"{u}}$nneth   formula \cite{SaW}. Let
$X$ be the closure of $Y$ and $D$ be the effective boundary 
divisor, then $\kappa(D, X)=1$ \cite{Zh1}. 
By \cite{Zh1}, $Y$ contains no complete curves and the boundary 
$X-Y$ is connected.  It is obvious that we have a surjective
morphism from $Y$ to $C$ (the projection). 
  But $Y$ is not affine since  $\kappa(D, X)=1$. 

\end{remark}
\begin{corollary} Let $Y$ be a  Stein algebraic surface. 
Let $X$ be its  completion. Then $Y$ is affine if and only if 
$\kappa(D, X)=2$, where $D$ is an effective divisor with support $X-Y$. 
Moreover, if $Y$ is smooth, then $Y$ is not affine if and only if 
$\kappa(D, X)=0$.
\end{corollary}

{\it Proof}.  We may assume that  $Y$ is normal as before.
 Since $Y$ is Stein,  it contains no complete
curves and  the boundary is connected and of pure codimension 1 \cite {N}.
Now the first claim  is obvious by Theorem 2.4. 

For the second claim, notice that if $\kappa(D, X)=1$,
then we have surjective morphism from $Y$ to a smooth affine
curve $C$. By Corollary 2.7,   $Y$ is affine which contradicts  the fact  
$\kappa(D, X)=1$. 
\begin{flushright}
 Q.E.D. 
\end{flushright}

\begin{remark}   If  $Y$ is a smooth algebraic 
Stein variety with dimension 3, then $\kappa(D, X)\neq 2$ 
(the following Theorem 3.4) but 
$\kappa(D, X)=1$ is possible \cite{Zh1}.  The case $\kappa(D, X)=0$
is a mystery. I do not know the existence of such algebraic variety.
In surface case, J.-P. Serre gave an example (\cite{H2},
Chapter VI, Example 3.2). The open surface
$Y$ is Stein but $\kappa(D, X)=0$ (\cite{Ku}, Lemma 1.8)
 
\end{remark}

\section{Threefolds}

We need the following theorems proved in \cite {Ku, Zh1, Zh2}. 

\begin{theorem}{\bf [Mohan Kumar]}   Let
$Y$ be a smooth algebraic surface over $\Bbb{C}$
with   $H^i(Y, \Omega^j_Y)=0$   for all $j\geq 0$ and $i>0$,
then  $Y$ is one of the following 
 
      (1) $Y$ is affine.

      (2) Let $C$ be an elliptic curve and $E$ the unique nonsplit 
      extension of $\mathcal{O}$$_C$ by itself.  
      Let ${X=\Bbb{P}}_C(E)$ and  $D$ be the canonical section, then $Y=X-D$.

       (3) Let $X$ be a projective rational surface with an effective 
       divisor $D=-K$ with $D^2=0$, $\mathcal{O}$$(D)|_D$ be nontorsion and 
       the dual graph of $D$ be $\tilde{D}_8$ or $\tilde{E}_8$, then $Y=X-D$.
\end{theorem}

For the above  type (2) and type (3) projective  surfaces, 
$\kappa (D, X)=0$ since 
$H^0(X, \mathcal{O}(nD))=\Bbb{C}$ for all nonnegative integer $n$
(\cite{Ku}, Lemma 1.8).

\begin{theorem}{\bf [Zhang]}  If  $H^i(Y, \Omega^j_Y)=0$ for all $j\geq 0$, 
     $i>0$ and 
$H^0(Y, \mathcal{O}$$_Y)\not=\Bbb{C}$, then    we have

        (1) There is a smooth projective 
curve $\bar{C}$,  and a smooth  affine curve $C$ such that 
the following diagram commutes

\[
  \begin{array}{ccc}
    Y                           &
     {\hookrightarrow} &
    X                                 \\
    \Big\downarrow\vcenter{%
        \rlap{$\scriptstyle{f|_Y}$}}              &  &
    \Big\downarrow\vcenter{%
       \rlap{$\scriptstyle{f}$}}      \\
C        & \hookrightarrow &
\bar{C}
\end{array}
\]
where f is proper and surjective, every fibre of f over $\bar{C}$
 is connected, 
 general fibre  is  smooth. Also general fibre 
 of $f|_Y$ is connected  and smooth.

       (2) Every fibre $S$ of $f|_Y$ over $C$  satisfies the same vanishing 
       condition, i.e., 
$H^i(S, \Omega^j_S)=0$.

       (3) We cannot have mixed types of  fibres, i.e., all  fibres must be 
exactly   one 
of the three  types of surfaces in Mohan Kumar's above classification. 

\end{theorem}

\begin{theorem} If $Y$ is a smooth threefold with 
 $H^i(Y, \Omega^j_Y)=0$
for all $j\geq 0$,  $i>0$ and $X$ is  a smooth completion
of $Y$  such that 
$X-Y$ is the support of 
    an effective boundary  divisor $D$
with  simple  normal crossings,
then  $\kappa (D, X)\neq 2$.
\end{theorem}
{\it Proof}.  
Suppose   $\kappa (D, X)  > 1$, then we have 
the above commutative diagram in Theorem 3.2. 
For a general point  $t\in C$,
 let  $X_t=f^{-1}(t)$  be the corresponding  smooth projective fibre
on $X$.
Let $D_t=D|_{X_t}$ be the restriction of $D$ on the fibre
$X_t$.  If  
$\kappa (D, X)>1$, then  by inequality (2.5),
we have 
$$  2\leq \kappa (D, X) \leq  \kappa (D_t, X_t)+1. 
$$
Thus $\kappa (D_t, X_t)\geq 1$.  Since 
$H^i(S_t, \Omega^j_{S_t})=0$  for all $j\geq 0$ and $i>0$ ,
$\kappa (D_t, X_t)=2$ by Theorem  3.1 and Lemma 1.8 \cite {Ku}. 
Pick  a  prime  divisor  $D_1 \leq  D$  such that
$D_1\cap f^{-1}(C)$ 
is empty, i.e.,  $f(D_1)=t_1\in \bar{C}-C$.
By the  Riemann-Roch  formula, for  the general point 
$t\in C$, 
$$  h^0(nt_1-t)-h^1(nt_1-t)={\mbox{deg}}(nt_1-t)+1-g(\bar{C})=n-g(\bar{C}),
$$ 
where  $g(\bar{C})$  is the genus of $\bar{C}$
and $h^i(nt_1-t)=h^i(\bar{C}, {\mathcal{O}}_{\bar{C}}(nt_1-t))$,
i=1, 2.
Choose  $n>g+1$, then  $h^0(nt_1-t)>1.$ Choose a nonconstant 
section $\sigma$ in $H^0(\bar{C}, {\mathcal{O}}_{\bar{C}}(nt_1-t))$, then
  $${\mbox{div}}\sigma +nt_1-t\geq 0.
$$ 
Pull it back to $X$, we have  
$$ f^*({\mbox{div}}\sigma +nt_1-t)
={\mbox{div}}(f^*\sigma)+nD_1-F\geq 0,
$$
where  $D_1=f^*(t_1)$  and   $F=f^*(t)$.  
So there is a nonzero section $f^*\sigma$ in $H^0(X, {\mathcal{O}}(nD_1-F))$.
Since  $H^0(X, {\mathcal{O}}(nD_1-F))\subset H^0(X, {\mathcal{O}}(nD-F))$,
by  Lemma 2.6,  
$$\kappa (D, X)=\kappa (D_t, X_t)+1=3.
$$ 
Therefore  $\kappa (D, X)\neq  2.$
\begin{flushright}
 Q.E.D. 
\end{flushright}

\begin{theorem} If $Y$ is a smooth algebraic Stein 
variety with dimension 3, then 
$\kappa(D, X)\neq 2$. 
\end{theorem}

$Proof.$  If $\kappa (D, X)> 1$, then we have  a 
surjective                    morphism from $Y$
to a smooth affine curve $C\subset \Bbb{C}$ defined  by 
a  nonconstant  regular   function   $f\in  H^0(Y, {\mathcal{O}}_Y)$.
It gives a rational map from $X$ to $\bar{C}$, the closure of 
$C$.
By Hironaka's   elimination   of indeterminacy and Stein factorization, we 
have the commutative diagram similar to Theorem 3.2
\[
  \begin{array}{ccc}
    Y                           &
     {\hookrightarrow} &
    X                                 \\
    \Big\downarrow\vcenter{%
        \rlap{$\scriptstyle{f|_Y}$}}              &  &
    \Big\downarrow\vcenter{%
       \rlap{$\scriptstyle{f}$}}      \\
C        & \hookrightarrow &
\bar{C},
\end{array}
\]
where $f$ is proper and surjective, every fibre $X_t=f^{-1}(t)$
 over $t\in\bar{C}$
 is connected, 
 general fibre  is  smooth. Also general fibre 
 of $f|_Y$ is connected  and smooth.  
Since $Y$  is Stein,  every  open  fibre  $S_t=X_t\cap  Y$
is Stein.  Since  $\kappa (D, X)> 1$, 
$$ 2\leq   \kappa (D, X) \leq  \kappa (D_t, X_t)+1,
$$
where $D_t=D|_{X_t}$ is the restriction divisor on the surface 
$X_t$.
Thus   $\kappa (D_t, X_t)=2$  by  Corollary 2.10. By Theorem 2.4,
$S_t$   is   affine.   By Lemma  2.6, 
$\kappa (D, X)=3$.
\begin{flushright}
 Q.E.D. 
\end{flushright}
\begin{theorem} {\bf[Zariski's Main Theorem]} 
Let $X$ be a normal variety  over a field $k$ and
let $f:X'\longrightarrow X$
be a birational morphism with finite 
fibre from a variety $X'$ to $X$.
Then   $f$  is an isomorphism of  $X'$
with an open subset $U\subset X$. 
\end{theorem}

\begin{theorem}{\bf [Neeman]} Let  $X=$${\mbox{Spec}}$A be a scheme,
$U\subset X$  a quasi-compact  Zariski  open  subset. Then
$U$  is affine  if and only if  $H^i(U, {\mathcal{O}}_U)=0$
for $i\geq 1$.
\end{theorem}

\begin{theorem} Let $Y$   be  a smooth threefold with 
 $H^i(Y, \Omega^j_Y)=0$
for all $j\geq 0$ and $i>0$.  Let  $X$   be   a smooth completion
of $Y$  such that 
  $D$ is  an effective boundary divisor on $X$
with  simple  normal crossings. 
If   $\kappa (D, X)>1$, then the associated scheme of  $Y$  is  isomorphic
to {\mbox{Spec}}$\Gamma(Y, {\mathcal{O}}_Y)$.  
\end{theorem}
{\it Proof}.  Since $\kappa (D, X)>1$,  the $D$-dimension
is 3  by Theorem 3.4. So we have a dominant morphism 
$h$  from $Y$  to
 $\Bbb{A}^3$ defined by three algebraically independent nonconstant functions
on $Y$. Let $Z$ be the normalization of $\Bbb{A}^3$ in $Y$, then we
have a morphism $g$ from  $Y$ (we identify $Y$ with its associated Scheme,
\cite{H1}, Chapter II, Proposition 2.6) to $Z=$Spec$\Gamma(Y, {\mathcal{O}}_Y)$
and  $g$ is birational since $Y$ and $Z$ have  the same function field. 
We claim that 
$g$ is injective. 
If there are two points $y_1$ and $y_2$ in $Y$ such that
they are mapped to the same point $z$ by $g$, i.e., $g(y_1)=g(y_2)=z$,
then for any regular function $r$ on $Y$, $r(y_1)=r(y_2)$.
Since $\kappa (D, X)\geq2$,
by Ueno's construction  (\cite{Uen}, page 46), there is a proper surjective morphism from $X'$
to ${\Bbb{P}}^1$, where $X'$ is a new smooth projective threefold 
obtained by resolving the indeterminacy, i.e., we have a surjective proper
birational morphism from $X'$ to $X$ such that it is an isomorphism on $Y$.
We still call the new threefold  $X$. So we have the 
following commutative  diagram
 \[
  \begin{array}{ccc}
    Y                           &
     {\hookrightarrow} &
    X                                 \\
    \Big\downarrow\vcenter{%
        \rlap{$\scriptstyle{f|_Y}$}}              &  &
    \Big\downarrow\vcenter{%
       \rlap{$\scriptstyle{f}$}}      \\
C        & \hookrightarrow &
{\Bbb{P}}^1
\end{array}
\]
where $f$ is proper and surjective, every fibre
$X_t=f^{-1}(t)$  on  $X$
 is connected, 
and general  fibre of $f$ over ${\Bbb{P}}^1$
   is  smooth and irreducible. By Theorem 3.2, every  open fibre on $Y$
satisfies the same vanishing condition.

  Since the $D$-dimension of $X$ is 3, 
by Equation (2.5), for a general 
fibre $X_t=f^{-1}(t)$, $\kappa (D_t, X_t)=2$, where $D_t=D|_{X_t}$.
Let $S_t=X_t-D_t$. 
 Since $H^i(S, \Omega^j_{S_t})=0$, the boundary divisor  $D_t=X_t-S_t$
is connected and $S_t$ contains no complete curves
\cite{Ku,Zh1,Zh2}.  These  three 
conditions:    
(1) $S_t$
has no complete curves;
(2) the boundary divisor $D_t$ is connected and 
(3)
$\kappa(D_t, X_t)=2$
imply 
 that $S_t$ is affine  (Theorem 2.4).

   We claim  that   if  $t_1$ and $t_2$ are two distinct points on $C$ and 
$y_1\in S_{t_1}=f^{-1}(t_1)\cap Y$, $y_2\in S_{t_2}=f^{-1}(t_2)\cap Y$, 
then there is  a regular function $R$
on $Y$, such that  $R(t_1)\neq R(t_2)$. In fact, we can take $R$ to be $f$.
     Now suppose that $y_1$ and $y_2$  are two distinct  points  on the 
same connected affine  surface $S_a$, where $a\in C$. 
Then  $f(y_1)=f(y_2)=a$. Since $S_a$ is affine, there 
is a regular function  $r$ on $S_a$ such that 
$r(y_1)=1$ but $r(y_2)=0$.
From the short exact sequence
$$ 0\longrightarrow 
 {\mathcal{O}}_Y
\longrightarrow 
 {\mathcal{O}}_Y
\longrightarrow 
 {\mathcal{O}}_{S_a}
\longrightarrow 
0,
$$
where the first map is defined by $f-a$, 
we
have surjective  map from  $H^0(Y, {\mathcal{O}}_Y)$ to 
$H^0(S, {\mathcal{O}}_S)$ since  $H^1(Y,  {\mathcal{O}}_Y)=0$.
So the regular function
$r$ on $S_a$ can be lifted to $Y$, i.e, there is a
regular function $R$ on $Y$, such that $R|_{S_a}=r$. Thus  $R(y_1)=1$
but $R(y_2)=0$. 
Therefore  $g$ must be injective on   every  fibre $S_t$, $t\in  C=f(Y)$.

     Since  $g$ is an injective birational morphism,
 by Zariski's Main
Theorem  (\cite{Mu1}, Chapter 3, Section 9), $g$ is an open immersion from  
$Y$  
to ${\mbox{Spec}}\Gamma(Y, {\mathcal{O}}_Y)$. By  Neeman's 
theorem  \cite{N}, $Y$  is affine. 
\begin{flushright}
 Q.E.D. 
\end{flushright}

\begin{theorem}Let $Y$   be  a smooth threefold with 
 $H^i(Y, \Omega^j_Y)=0$
for all $j\geq 0$ and $i>0$.  Let  $X$   be   a smooth completion
of $Y$  such that 
  $D$ is  an effective boundary divisor on $X$
with  simple  normal crossings.
Then $Y$ is affine  if and only if    $Y$ is 
regularly separable.
\end{theorem}
{\it Proof}. If $Y$ is affine, then   $Y$
is a closed subset of  ${\Bbb{C}}^n$
for some  $n\in  \Bbb{N}$. 
There is a polynomial   $f$ 
on  ${\Bbb{C}}^n$  such  that
for two  distinct  points  $y_1$ and $y_2$  in $Y$,
$f(y_1)\neq  f(y_2)$. Obviously $f$ is a regular function 
on $Y$.
 So $Y$  is regularly separable.

  Suppose now that 
$Y$ is 
regularly separable.   
Since $D$ is effective, $h^0(X, {\mathcal{O}}_X(nD))>0$
for all $n\geq 0$. 
Since $Y$ is 
regularly separable, the $D$-dimension of $X$ is at least 1.
By the proof of  Theorem 3.4, we have the same fibre space and commutative 
diagram\[
  \begin{array}{ccc}
    Y                           &
     {\hookrightarrow} &
    X                                 \\
    \Big\downarrow\vcenter{%
        \rlap{$\scriptstyle{f|_Y}$}}              &  &
    \Big\downarrow\vcenter{%
       \rlap{$\scriptstyle{f}$}}      \\
C        & \hookrightarrow &
\bar{C}.
\end{array}
\]
 For a  general  point $t$ in $C$ such that 
the fibre $X_t=f^{-1}(t)$ is smooth and irreducible,
we know that there are 3 possible surfaces as in Theorem 3.1.
Since $Y$ is regularly separable, the open   surface  $S_t=X_t|_Y$  
  is   affine (\cite{Ku}, Lemma 1.8). Let $D_t=D|_{X_t}$ be the 
boundary divisor supported in $X_t-S_t$,
then  the $D_t$-dimension is 2. By upper semi-continuity theorem,
for all points $t$ in $C$, the $D_t$-dimension is 2. By 
 Lemma 2.6, the $D$-dimension of  $X$ is 3. 
So   as  in the proof of Theorem 3.7, 
 we have  an  injective  birational  morphism $g$ from $Y$
to Spec$\Gamma(Y, {\mathcal{O}}_Y)$.
Thus   $Y$ is an open subset of an affine variety by Zariski's  Main Theorem.
By Neeman's Theorem, $Y$ is affine since 
$H^i(Y, {\mathcal{O}}_Y)=0$.
\begin{flushright}
 Q.E.D. 
\end{flushright}

\begin{remark} $H^i(Y, \Omega^j_Y)=0$
for all $j\geq 0$ and $i>0$ is a  necessary condition. 
 Theorem 3.7 does not hold if we drop this assumption. 
That is, the following four conditions 
are necessary for affineness but not sufficient

(1) $Y$ contains no complete curves;

 (2) For a smooth completion $X$ of $Y$, $X-Y$
 is connected and is of pure codimension 1;

(3) Let $D$ be the boundary divisor, the $D$-dimension of $X$ is 3;

(4) $Y$ is regularly separable. 

For the counter-example, see Remark 2.5.
\end{remark}

\begin{corollary} If $H^i(Y, \Omega^j_Y)=0$
for all $j\geq 0$  and $i>0$,
$H^0(Y, \mathcal{O}$$_Y)\not=\Bbb{C}$,  and $Y$ is not 
affine, then 
$\kappa(X)=-\infty$ and $\kappa(D, X)=1$
for  every smooth  completion $X$ of $Y$ 
and any  effective  boundary divisor $D$
with support  $X-Y$.  
\end{corollary}
{\it Proof}.  Since $Y$  is not affine,  
$\kappa(D, X)<3$  by Theorem 3.7. Then the claim is an immediate 
consequence of  
Theorem 3.3  and  Theorem 7 \cite{Zh1}. 
\begin{flushright}
 Q.E.D. 
\end{flushright}

\begin{theorem} An algebraic manifold  $Y$   is    affine   if and only if 
 $H^i(Y, \Omega^j_Y)=0$
for all $j\geq 0$, $i>0$ and  $Y$ is 
regularly separable. 
 \end{theorem}
{\it Proof}.  We will prove the theorem by 
 induction on the dimension. By Theorem 3.1 and 3.8, 
we may assume 
that  the theorem holds for dimension $d$ algebraic manifolds. Suppose  now
dim$Y=d+1$. Pick a non-constant regular function on $Y$, we can construct 
 a similar 
fibre space as in  Theorem 3.8, that is,  all fibres in $Y$ have dimension 
$d$ and
satisfy the same vanishing condition \cite{Zh1}. By inductive assumption,
since all fibres are  regularly separable,
every  connected     open fibre in $Y$   is   affine. By  Lemma 2.6, 
$\kappa(D, X)=d+1$. Now we can apply the proof of Theorem 3.8
to argue that the map $g$ is a birational injective morphism from $Y$ to 
Spec$\Gamma(Y, {\mathcal{O}}_Y)$ since  $\kappa(D, X)=d+1$
and $Y$ is 
regularly separable.  By Zariski's  Main Theorem, 
$Y$ is  isomorphic  to an open  subset of 
Spec$\Gamma(Y, {\mathcal{O}}_Y)$. 
  By Neeman's theorem, $Y$ is affine.
\begin{flushright}
 Q.E.D. 
\end{flushright}

\end{document}